\documentclass[12pt]{article}
\usepackage{amssymb}
\usepackage{latexsym}
\usepackage{cite}
\usepackage{graphicx}

\input{amssym.def}
\input{amssym.tex}
\def\C{{\mathbb C}}

\def\Z{{\mathbb Z}}
\def\mH{\mathop{\otimes}_{H}}
\def\D {/\!\!\!\!D}
\newcommand{\CC}{\hbox{{$\cal C$}}}

\newcommand{\extd}{{\rm d}}
\newcommand{\del}{\partial}

\newcommand{\tens}{\mathop{\otimes}}

\newcommand{\Ad}{{\rm Ad}}
\newcommand{\id}{{\rm id}}

\renewcommand{\>}{\rangle}

\newcommand{\proof}{{\bf Proof\ }}
\newcommand{\eproof}{$\quad \diamond$\bigskip}

\newcommand{\eqn}[2]{\begin{equation}#2\label{#1}\end{equation}}
\newcommand{\lcross}{{>\!\!\!\triangleleft}}

\newtheorem{lemma}{Lemma}[section]
\newtheorem{prop}[lemma]{Proposition}

\newtheorem{theo}[lemma]{Theorem}

\newtheorem{defi}[lemma]{Definition}

\begin{document}
\begin{center}
{NONCOMMUTATIVE RIEMANNIAN GEOMETRY \\OF THE ALTERNATING GROUP
${\cal A}_{4}$.}
\end{center}
\begin{center}
{\bf F. Ngakeu }$^{a,b,}$\footnote{supported by CIUF-CUD}, 
{\bf S. Majid $^{c,}$}\footnote{Royal Society University Research Fellow}
{\bf and  D. Lambert$^{d}$ }
\end{center}
\begin{center}
${\bf a}:$ Unit\'e de Physique Th\'eorique et de Physique 
Math\'ematique, U.C.L\\
 2, Chemin du Cyclotron, B-1348  Louvain-La-Neuve,
Belgium\\

${\bf b}:$ Institut de Math\'ematiques et de Sciences Physiques,\\
 BP 613 Porto-Novo, Benin. e-mail: fngakeu@yahoo.fr
 \bigskip \\
${\bf c}:$ School of Mathematical
Sciences, Queen Mary University of
London,\\ Mile End Rd, London E1 4NS, UK
 \bigskip \\
${\bf d}:$ Facult\'es des Sciences, F.U.N.D.P\\
 61 Rue de Bruxelles, B-5000
Namur, Belgium 
\end{center}
\begin{center}
July, 2001 
\end{center}

{\bf Abstract}: We study the noncommutative Riemannian geometry of
the alternating group $A_4=(\Z_2 \times \Z_2)\lcross\Z_3$ using
the recent formulation for finite groups in \cite{maj}. We find a
unique `Levi-Civita' connection for the invariant metric, and find
that it has Ricci-flat but nonzero Riemann curvature. We show that
it is the unique Ricci-flat connection on $A_4$ with the standard
framing (we solve the vacuum Einstein's equation). We also propose
a natural Dirac operator for the associated spin connection and solve
the Dirac equation. Some
of our results hold for any finite group equipped with a cyclic
conjugacy class of 4 elements. In this case the exterior algebra
$\Omega(A_4)$ has dimensions $1:4:8:11:12:12:11:8:4:1$ with top-form
9-dimensional. We also find the noncommutative cohomology
$H^1(A_4)=\C$.

\section{Introduction}\baselineskip 20pt
\qquad  A constructive formalism of noncommutative Riemannian
geometry has recently been developed in \cite{maj-br} and
\cite{maj}, using quantum group methods. Here a general and
possibly noncommutative algebra or `coordinate ring' is equipped
with a `quantum Riemannian manifold' structure consisting of a
frame bundle (with quantum group fiber) to which the differential
calculus\cite{w} or exterior algebra bundle is associated. The
approach builds on the established  theory of quantum principal
bundles\cite{tbmaj} and adds to this notions of `framing',
`coframing' (or metric) and Levi-Civita type metric-compatible
connections with Riemann and Ricci curvature.

This approach applies also to finite sets and finite groups, and
allows one to endow them with nontrivial Riemannian manifold
structures. This would not be possible within conventional
differential geometry, but noncommutative differential structures
are more general even when the coordinate ring is commutative, and
allow a rich structure even for finite sets. Particularly, for a
finite group there are natural choices for translation biinvariant
differential structure, namely labelled by the nontrivial
conjugacy classes (see \cite{w}). There is a natural frame bundle,
namely with fiber another copy of the group\cite{maj-br}, and
there is a natural `Killing form' inducing a metric\cite{maj},
which may, however, be degenerate.  So far, only the case of $S_3$
in \cite{maj} has been worked out in detail it is shown there that
this has the natural structure of a noncommutative Einstein
manifold with unique Levi-Civita connection for the Killing
metric, and with  Ricci curvature essentially proportional to the
metric. Some other aspects of the noncommutative geometry of
 $S_3$ are in \cite{mr}.

In this paper we extend the repertoire of examples with a detailed
study of the alernating group $A_4$ from this point of view. It
turns out that this example is of key interest because, like
$S_3$, it has a natural invariant metric with unique Levi-Civita
connection, but this connection is Ricci-flat. Thus, it provides
the first concrete example of the solution of the vacuum Einstein
equations in this theory. The model is also interesting because
the group is nonabelian enough to have nontrivial curvature (the
Riemann tensor does not vanish) but simple enough to be fully
computable.

Some of the computations are done without reference to the actual
group and apply to any group with similar `cyclic' conjugacy
class. We include for example $SL_2(\Z_3)$ in this family. After
some preliminaries, we start in Section~3 at this general level.
We find the Woronowicz exterior algebra associated to the
conjugacy class and show that it is {\em not} in fact
quadratically generated. It is in fact a good example where we see
the absolute necessity of relations in higher degree. We then find
the unique form of the invariant metric, namely
\[ \eta^{a,b}=\delta_{a,b}+\mu\]
in a suitable basis, with $\mu$ a parameter. We also characterise
the torsion free and cotorsion free regular connections. In
Section~4 we specialise to $A_4$ and find the associated
Levi-Civita or metric-compatible connection and its Ricci
flatness. We also show that it is the unique regular Ricci-flat
connection on $A_4$ independently of metric compatibility.

In Section~5 we look at the Dirac operator appropriate to the
metric, as a step towards comparison with Connes' approach to
noncommutative geometry\cite{con}. As in \cite{maj} for $S_3$, it
is not naturally hermitian but is fully diagonilisable. Finally,
Section 6 contains some further results including the
noncommutative de Rahm cohomology and a link to the differential
(but not Riemannian) geometry of $S_4$ in \cite{maj-per}, as well
as concluding remarks.

We note that following \cite{maj-br} there have been some other
attempts at a Riemannian geometry on finite groups, but different
from the one used here, see for instance \cite{dim, kbres}  and
references therein. While the use of a conjugacy class to define
an exterior algebra, i.e. the notion of differential structure, is
the same in all approaches following \cite{w}, see for instance
\cite{sitarz}, the formulation of \cite{maj-br}\cite{maj} is the
only one that features some kind of metric compatibility or
`Levi-Civita' notion for the spin connection, as well as the only
one that features a nontrivial (and nonuniversal) differential
structure in the fiber direction, and hence an actual
noncommutative geometry of the frame bundle. Other problems solved
in \cite{maj} were the formulation and computation of the Ricci
tensor, see \cite{maj-con} for a discussion.

\section{Preliminaries}

\qquad In this section, we briefly recall the basic definitions of
differential structures\cite{w} and of noncommutative Riemannian
geometry\cite{maj}, specialised to the case of  finite groups that
we need in the paper. Thus, we work with the Hopf algebra
$H=\C[G]$ of functions on a finite group $G$. We equip it with its
standard basis $(\delta_{g})_{g\in G}$ defined by Kronecker
delta-functions $\delta_{g}(h)=\delta_{g}^{h}, \forall g,h \in G$.
Let ${\cal C}$ be a nontrivial conjugacy class of $G$, and
$\Omega_{0}$ the vector subspace
\begin{eqnarray}
\Omega_{0}=\{ \delta_{a}|\ a\in {\cal C}\}=\C {\cal C}.
\end{eqnarray}
Any $\Ad$-stable set not containing the group identity can be used
here, but we focus on the irreducible case of a single conjugacy
class. {}From the Woronowicz's theory in \cite{w}, the first order
differential calculus associated to ${\cal C}$ is generated over
$\C[G]$ by $\Omega_0$ and given by
\begin{eqnarray}\label{extdif}
\extd f=\sum\limits_{a \in {\cal C}}(\partial^{a}f)e_{a},\quad
 \partial^{a}=R_{a}-\id,\quad e_{a}f=R_{a}(f)e_{a}
\end{eqnarray}
$\forall a \in {\cal C}, f \in H $, where the operator $R_{a}$ is
defined by $R_{a}(f)(g)=f(ga), \forall g \in G$. The Maurer-Cartan
1-form $e:\Omega_{0} \longrightarrow \Omega^{1}(H)$ is given by
\begin{eqnarray}
e_{a}=e(\delta_{a})=\sum\limits_{g\in
  G}\delta_{g}\extd
\delta_{ga}
\end{eqnarray}
and higher differential forms are obtained from Woronowicz
skew-symmetrization procedure \cite{w}, using the braiding
\begin{eqnarray}
\Psi (e_{a} \mH e_{b})=e_{a b a^{-1}}\mH e_{a}
\end{eqnarray}
The Maurer-Cartan equation takes the form
\begin{eqnarray}\label{MCequation}
\extd e_{a}=\theta\wedge e_a+e_a\wedge \theta,\quad
\theta=\sum\limits_{a \in {\cal C}}e_a.
\end{eqnarray}
The element $\theta$ obeys $\theta\wedge\theta=0$ and generates
$\extd$ in general as graded-commutator with $\theta$. Lemma 5.3
in \cite{maj} gives the full set of relations of $\Omega^2(H)$,
namely
\begin{eqnarray}\label{lemdemaj}
\sum\limits_{a,b \in {\cal C};\;a b =g}\lambda _{a}^{g,\beta}
e_{a}\wedge e_{b}=0, \forall g \in G, \forall \beta.
\end{eqnarray}
where for $g\in G$ fixed, $\{\lambda^\beta\}$ is a basis of the
invariant subspace of the vector space spanned by $\CC\cap
g\CC^{-1}$ under the automorphism $\sigma(a) =a^{-1}g$. There are
also cubic and higher degree relations (which are in fact
nontrivial in our case of $A_4$) but we will not need them
explicitly (most of Riemannian geometry needs only 1-forms and
2-forms.)

Next, following \cite{maj}, a framing means a basis of $\Omega^1(H)$
over $\C[G]$, and an action of the frame group. In our case we
chose the framing to be the components $\{e_a\}$ of the
Maurer-Cartan form as above and for frame group we choose $G$
itself, acting by $\Ad$. This is a canonical choice and its
classical meaning is explained in \cite{maj-br}. A spin connection
is then a collection $\{A_{a}\}_{a \in {\cal C}}$ of component
1-forms. Its associated covariant derivative is defined on an
1-form $\alpha= \alpha^{a}e_{a}$ by
\begin{eqnarray}\label{covariantderivative}
\nabla \alpha= \extd \alpha^{a} \mH e_{a}-\alpha^{a}\sum\limits_{b
  \in {\cal C}}A_{b}\mH (e_{b^{-1}a b}-e_{a})
\end{eqnarray}
with summation on $a$. The associated torsion tensor $T:\Omega
^{1}(H)\longrightarrow \Omega^{2}(H)$ is defined by $T\alpha
=\extd\wedge \alpha -\nabla \alpha$ and the zero- torsion
condition is vanishing of
\begin{eqnarray}\label{torsionequation}
\bar{D}_{A}e_{a}\equiv \extd e_{a}+\sum\limits_{b
  \in {\cal C}}A_{b}\wedge (e_{b^{-1}a b}-e_{a}),\quad \forall
  a \in {\cal C}
\end{eqnarray}
The spin connection here has values in the dual space
$\Omega_0^*$, which is a `braided-Lie algebra' in  a precise
sense. Associated to this geometrical point of view, there is a
regularity condition
\begin{eqnarray}\label{regularitycondition}
\sum\limits_{a, b \in {\cal C};\;a b =g}A_{a}\wedge A_{b}=0 ,\quad
\forall g\neq e,g\notin {\cal C}.
\end{eqnarray}
The curvature $\nabla^{2}$ associated to a regular connection $A$
is in frame bundles terms a collection of 2-forms $\{F_{a}\}_{a\in
{\cal C}}$ defined by
\begin{eqnarray}\label{2formdecourbure}
F_{a}=\extd A_{a}+\sum\limits_{c,d \in {\cal C}, cd=a}A_{c}\wedge
A_{d}- \sum\limits_{c\in {\cal C}}(A_{c}\wedge A_{a}+A_{a}\wedge
A_{c})
\end{eqnarray}
while the Riemann curvature ${\cal R}: \Omega^{1}(H)
\longrightarrow \Omega^{2}(H)\mH \Omega^{1}(H)$ is given by
\begin{eqnarray}\label{formuledelacourbureriem}
{\cal R}\alpha =\alpha^{a}\sum\limits_{b \in {\cal C}}F_{b} \mH
(e_{b^{-1}a b}-e_{a})
\end{eqnarray}
Finally, the Ricci tensor is given by
\begin{eqnarray}\label{ricci}
Ricci=\sum\limits_{a,b,c \in {\cal C}}i(F_c)^{a
  b}e_{b}
\mH (e_{c^{-1}a c}-e_{a})
\end{eqnarray}
where $i(F_{c})=i(F_c)^{a
  b}e_{a}\mH e_{b}$ and $i:\Omega^{2}(H) \longrightarrow
\Omega^{1}(H)\mH \Omega^{1}(H)$ is a lifting which splits the
projection of $\wedge $. A canonical choice is \cite{maj}
\begin{eqnarray}\label{canonicallift}
i(e_{a}\wedge e_{b})=e_{a} \mH e_{b}-\sum\limits_{\beta}
\gamma^{\beta ,a}\sum\limits_{c,d \in {\cal
    C},\;cd=ab}\lambda_{c}^{\beta}
e_{c}\mH e_{d}
\end{eqnarray}
where $\{\gamma^{\beta\}}$ are the dual basis to the
$\{\lambda^{\beta}\}$ with respect to the dot product as vectors
in $\C {\cal C}\cap ab{\cal C}^{-1}$. Another canonical `lift' is
$i'=\id-\Psi$ but note that in this case $i'\circ\wedge$ is {\em
not} a projection operator.

There are two further structures that one may impose in this
situation. First of all, given a choice of framing $\{e_a\}$, a
metric $g$ is defined as a coframing $\{e^{*a}\}$, i.e. again a
basis of $\Omega^1$ but now as a right $\C[G]$-module. and
transforming under the contragradient action of $G$.(The
corresponding metric is $g=\sum_a e^{*a}\tens_H e_a$). The
cotorsion of a spin connection is the torsion with respect to the
coframing, and is given by
\begin{eqnarray}\label{cotorsion}
D_{A}e^{*a}\equiv \extd e^{*a}+\sum\limits_{b \in {\cal C}} (e^{*
b a b^{-1}}-e^{*a})\wedge A_{b}.
\end{eqnarray}
Vanishing of cotorsion has the classical meaning of a
generalisation of metric compatibility of the spin connection,
see\cite{maj-br}. So we are usually interested in regular
torsion-free and cotorsion-free connections.

Finally, a `gamma-matrix' is defined\cite{maj} as a collection of
endomorphisms $\{\gamma _{a}\}_{a \in {\cal C}}$ of a vector space
$W$ on which $G$ acts by a representation $\rho_W$, a `spinor
field' is a $W$-valued function on $G$ and the Dirac operator on
the spinor fields is
\begin{eqnarray}\label{operateurDirac}
\D =\partial^{a}\gamma_{a}-A_{a}^{b}\gamma_{b}\tau_{W}^{a},
\end{eqnarray}
where $A_{b}=A_{b}^{a}e_{a}$ and $\tau_{W}^{a}=\rho_W(a^{-1}-e)$.
There is a canonical choice where $\gamma$ is built from $\rho_W$
itself, explained in \cite{maj}.

\section{Cyclic Riemannian structures}

In this section, we construct  Riemannian geometry on groups
endowed with conjugacy classes which obey a certain cyclicity
condition. For the case when the differential calculus is of
degree four, we determine the entire exterior algebra and the
moduli space of torsion free connections, and for any degree
$n\geq 2$, we give the general form of the invariant metric.

\begin{defi}\label{cyclicclass}
Let ${\cal C}$ be a conjugacy class with $n$ elements, $n\geq2$,
in a group $G$. We say that ${\cal C}$ is `cyclic' if there exists
at least one $t$ in  ${\cal C}$ such that $\Ad_{t}$ is a cyclic
permutation of ${\cal C}-\{t\}$ and the map $a\rightarrow
\Ad_{a}(t)$ is a permutation of ${\cal C}$.
\end{defi}

For $n=4$ we have the following characterisation of
$\Omega^{2}(H)$

\begin{prop}\label{prop3.1}
For a cyclic conjugacy class ${\cal C}=\{t,x,y,z\}$ of order 4 in
a finite group $G$, the bimodule $\Omega^{2}(H)$ of 2-form is
8-dimensional over $\C[G]$ and is defined by the following
equations
\begin{eqnarray}\label{2-form1}
e_{a}\wedge e_{a}=0, \quad \sum\limits_{a,b\in {\cal C};\;
ab=g}e_{a}\wedge e_{b}=0
\end{eqnarray}
$\forall a \in {\cal C}, g \in G$, where $(e_{a})_{a\in {\cal C}}$
is the basis of Maurer-Cartan 1-forms.
\end{prop}
\proof :\\
We assume the existence of an element $t\in\CC$ as in
Definition~\ref{cyclicclass}. Without loss of generality, we
denote the other elements of ${\cal C}$ by $x,y,z$ with the
following table for $\Ad$:
\begin{eqnarray}\label{tableau1}
\begin{tabular}{|c|c|c|c|c|}
\hline
$\Ad$ &t&x&y&z\\
\hline
$t$ & t & z & x & y \\
\hline
$x$&y&x&z&t\\
\hline
$y$&z&t&y&x\\
\hline
$z$ &x&y&t&z\\
\hline
\end{tabular}
\\ \nonumber
 {\mbox{Table 1}}
\end{eqnarray}
It follows that
\begin{eqnarray}\label{produits1}
tx=zt=xz;\quad ty=xt=yx;\quad tz=yt=zy;\quad xy=zx=yz
\end{eqnarray}
Using relations (\ref{produits1}), we apply the Woronowicz
antisymmetrization procedure to obtain the following relations of
the form (\ref{2-form1}) in $\Omega^{2}(H)$:
\begin{eqnarray}\label{les2-form}
e_{a}\wedge e_{a}&=&0,\quad \forall a \in {\cal C} \\ \nonumber
e_{t}\wedge e_{x}&+&e_{x}\wedge e_{z}+e_{z}\wedge e_{t}=0 \\
\nonumber e_{x}\wedge e_{t}&+&e_{t}\wedge e_{y}+e_{y}\wedge
e_{x}=0 \\ \nonumber  e_{t}\wedge e_{z}&+&e_{z}\wedge
e_{y}+e_{y}\wedge e_{t}=0 \\ \nonumber e_{x}\wedge
e_{y}&+&e_{y}\wedge e_{z}+e_{z}\wedge e_{x}=0 \\ \nonumber
\end{eqnarray}
This form (\ref{2-form1}) holds for any group since the elements
$e_a\tens e_a$ and $\sum_{ab=g}e_a\tens e_b$ are in the kernel of
$\id-\Psi$. However, using (\ref{lemdemaj}) one may see that they
are the {\em only} relations of $\Omega^{2}(H)$ which is therefore
of dimension 8 as stated. \eproof

{}From now, we choose a basis of $\Omega^{2}(H)$ to be
\begin{eqnarray}\label{base2-form}
\{e_{t}\wedge e_{x},\; e_{t}\wedge e_{y},\; e_{t}\wedge e_{z},\;
e_{x}\wedge e_{t},\; e_{y}\wedge e_{t},\; e_{x}\wedge e_{y},\;
e_{y}\wedge e_{z},\; e_{x}\wedge e_{z}\},
\end{eqnarray}
For convenience, we will sometime use indexes 1, 2, 3, 4 to refer
to $t, x, y, z$ respectively.

\begin{prop} In the setting of Proposition~\ref{prop3.1}, the
dimensions of the Woronowicz exterior algebra $\Omega(H)$ are
$1:4:8:11:12:12:11:8:4:1$ with top-form of degree 9. This algebra is
{\em not} quadratic, having additional relations in degree $\ge
6$.
\end{prop} \proof:\\
As above, we do not need the group itself but only the matrix for $\Ad$ 
restricted to the conjugacy class (i.e. Table~1). In fact we are computing the 
invariant part $\Lambda$ of the exterior algebra, with 
$\Omega(H)=H\tens\Lambda$ as a vector space.
This $\Lambda$ is generated over $\C$ by the $\{e_a\}$ with relations 
determined by
the braiding $\Psi$. Namely we set to zero the kernels of the antisymmetizers 
$A_m$ for $m\ge 2$. These $A_m$ are described in \cite{w} as a signed sum over 
permutations of $\{1,\cdots, m\}$ with transposition replaced by $\Psi$. 
This is not very convenient for computation and we employ instead a different 
but equivalent definition of the $A_m$ coming out of the theory of braided 
groups\cite{maj-book}. As recently discussed in \cite{maj-per}, we use the 
braided-integers 
\[ [m,-\Psi]=\id-\Psi_{12}+\Psi_{12}\Psi_{23}-\cdots\pm \Psi_{12}
\cdots\Psi_{m-1,m}=\id-\Psi_{12}(\id\tens [m-1,-\Psi]),\]
where $\Psi_{12}$ denotes $\Psi$ acting in the first and second places of 
$\Omega_0^{\tens m}$, etc. Then 
\[ A_m=[m,-\Psi]!=(\id\tens [2,-\Psi])(\id\tens [3,-\Psi])
\cdots [m,-\Psi].\]
In the braided groups approach to the exterior algebra we set to 
zero the kernels
of all these braided factorials.  It is straightforward to program 
these inductive
definitions. We first compute the $16\times 16$ matrices $\Psi$ 
acting in the 
tensor product basis $e_a\tens e_b$ and then the $A_m$ as above, 
up to $A_6$. 
The dimensions $\Omega^m(H)$ over $H$ are then $4^m-\dim\ker A_m$ 
and found to be as stated.
{}From the general form expected for the exterior algebra we assume
 the remaining dimensions 
for $A_7,A_8,A_9$ without explicit computation. Finally, the 
quadratic exterior algebra
is defined by setting to zero only the kernel of $A_2=\id-\Psi$ 
without additional
relations in higher degree. In that case in degree $m$ we set to 
zero the union 
of the null spaces $\id-\Psi_{12},\cdots,\id-\Psi_{m-1,m}$. Here we 
find dimensions
$1:4:8:11:12:12:12:\cdots$ i.e., fewer relations in degree $\ge 6$ 
(it appears that the
quadratic one is in fact infinite-dimensional). \eproof

In fact for most geometric purposes we need only the exterior
algebra up to degree 2, so we limit ourselves to the general
result about dimensions. In principle one may go on to compute 
explicit relations in higher degree and a Hodge * operator as 
in \cite{mr} using the metric below, etc. The
result is an important reminder that the degree 2 relations alone
may not be enough for a geometrically reasonable exterior algebra.

\begin{prop}\label{prop3.2}
\quad In the setting of Proposition \ref{prop3.1} above and for
the framing defined by the Maurer-Cartan 1-form, the moduli space
of torsion free connections is $3|G|$- dimensional and is given
by the following components 1-forms:
\begin{eqnarray}
\label{module1}
A_{t}&=&(1+\alpha)e_{t} +\gamma e_{x} +\lambda e_{y}+\beta e_{z}\\
\nonumber
A_{x}&=&\lambda e_{t}+(1+\beta)e_{x}+\alpha e_{y}+\gamma e_{z}\\
\nonumber
A_{y}&=&\beta e_{t}+\lambda e_{x}+(1+\gamma) e_{y}+\alpha e_{z}\\
\nonumber A_{z}&=& \gamma e_{t}+\alpha e_{x}+ \beta
e_{y}+(1+\lambda)e_{z}
\end{eqnarray}
where $\alpha,\; \beta,\; \gamma,\; \lambda$ are functions on $G$
such that
\begin{eqnarray}
\alpha+ \beta+ \gamma+\lambda=-1
\end{eqnarray}
Thus we have also
\begin{eqnarray}\label{sommedesAi}
\sum\limits_{a\in {\cal C}}A_{a}=0
\end{eqnarray}
\end{prop}
\proof:\\
We follow the same method as for $S_3$ in \cite{maj}. In the
framing defined by the Maurer-Cartan 1-form, the torsion free
connections obey the following equation (see eq.
(\ref{torsionequation}))
\begin{eqnarray}\label{eqtorsion}
\sum\limits_{b\in {\cal C}} A_{b}\wedge(e_{b^{-1}ab}-e_{a})+
\sum\limits_{b\in {\cal C}}(e_{b}\wedge e_{a}+e_{a}\wedge e_{b})=0
\end{eqnarray}
$\forall a \in {\cal C}$.
 Using Table~1, we write (\ref{eqtorsion}) as
\begin{eqnarray}\label{syseqtorsion}
A_{x}\wedge (e_{z}-e_{t})+ A_{y}\wedge (e_{x}-e_{t})+A_{z}\wedge
(e_{y}-e_{t})&{}&\\ \nonumber +(e_{x}+e_{y}+e_{z})\wedge
e_{t}+e_{t}\wedge (e_{x}+e_{y}+e_{z})&=&0
\\ \nonumber
A_{t}\wedge (e_{y}-e_{x})+ A_{y}\wedge (e_{z}-e_{x})+A_{z}\wedge
(e_{t}-e_{x})\\ \nonumber +(e_{t}+e_{y}+e_{z})\wedge
e_{x}+e_{x}\wedge (e_{t}+e_{y}+e_{z})&=&0\\\nonumber A_{t}\wedge
(e_{z}-e_{y})+ A_{x}\wedge (e_{t}-e_{y})+A_{z}\wedge
(e_{x}-e_{y})\\ \nonumber +(e_{t}+e_{x}+e_{z})\wedge
e_{y}+e_{y}\wedge (e_{t}+e_{x}+e_{z})&=&0\\\nonumber A_{t}\wedge
(e_{x}-e_{z})+ A_{x}\wedge (e_{y}-e_{z})+A_{y}\wedge
(e_{t}-e_{z})\\ \nonumber +(e_{t}+e_{x}+e_{y})\wedge
e_{z}+e_{z}\wedge (e_{t}+e_{x}+e_{y})&=&0\\\nonumber
\end{eqnarray}
 We just have to solve the first three equations since
the fourth one in this system can be obtained from the other by
simple summation. We set $A_{a}=A_{a}^{b}e_{b}$ (sum over
$b\in\CC$) for functions $ A_{a}^{b}\in H$ with
\begin{eqnarray}\label{potentiel}
A_{t}^{t}=1+\alpha\quad A_{x}^{x}=1+\beta,\quad
A_{y}^{y}=1+\gamma, \quad A_{z}^{z}=1+\lambda
\end{eqnarray}
We put this into the equations to be solved and write them in the
basis (\ref{base2-form}). Using the fact that each coefficient of
the basis element has to vanish, we obtain
\begin{eqnarray}\label{solpotentiel}
A_{x}^{t}&=&\lambda=A_{y}^{x},\quad  A_{y}^{t}=\beta=
A_{z}^{y}=A_{t}^{z}, \quad A_{z}^{t}=\gamma= A_{t}^{x}=A_{x}^{z}\\
\nonumber
 A_{z}^{x}&=&-1-\lambda-\gamma-\beta=A_{x}^{y}=A_{y}^{z},\quad
A_{t}^{y}=-1-\alpha-\beta-\gamma \\ \nonumber
&{}&\alpha+\beta+\gamma+\lambda=-1
\end{eqnarray}
as stated. Finally using these solutions one checks  by simple
computation that $A_{t}+A_{x}+A_{y}+A_{z}=0$.
\eproof \\

We now study the regularity of connections:
\begin{prop}
Under the hypothesis of Proposition~\ref{prop3.1} the regular
connections are either  solutions of the system:
\begin{eqnarray}\label{regularcase2}
A_{t}\wedge A_{t}&+& A_{x}\wedge A_{y}+A_{y}\wedge
A_{z}+A_{z}\wedge A_{x}=0\\ \nonumber A_{t}\wedge A_{x}&+&
A_{x}\wedge A_{z}+A_{y}\wedge A_{y}+A_{z}\wedge A_{t}=0\\
\nonumber A_{t}\wedge A_{y}&+& A_{x}\wedge A_{t}+A_{y}\wedge
A_{x}+A_{z}\wedge A_{z}=0\\ \nonumber A_{t}\wedge A_{z}&+&
A_{x}\wedge A_{x}+A_{y}\wedge A_{t}+A_{z}\wedge A_{y}=0
\end{eqnarray}
{\em or} solutions  of the system:
\begin{eqnarray}\label{regularcase1}
A_{t}\wedge A_{t}&=&0,\quad A_{x}\wedge A_{x}=0,\quad A_{y}\wedge
A_{y}=0,\quad A_{z}\wedge A_{z}=0\\ \nonumber A_{t}\wedge
A_{x}&+&A_{x}\wedge A_{z}+A_{z}\wedge A_{t}=0\\ \nonumber
A_{t}\wedge A_{y}&+&A_{x}\wedge A_{t}+A_{y}\wedge A_{x}=0\\
\nonumber A_{t}\wedge A_{z}&+&A_{y}\wedge A_{t}+A_{z}\wedge
A_{y}=0\\ \nonumber A_{x}\wedge A_{y}&+&A_{y}\wedge
A_{z}+A_{z}\wedge A_{x}=0
\end{eqnarray}
\end{prop}
\proof:\\
The general form of the regularity's equation is given by
(\ref{regularitycondition}). One then needs the multiplication
table at least for the elements of the class ${\cal C}$, by
enumeration of the cases we find that under the hypothesis of
Proposition~\ref{prop3.1}, the only possible cases are those shown
in Tables 2, 3. These correspond to the two possibilities
 stated. \eproof

\vspace{1cm} \hspace{1cm}
\begin{tabular}{|c|c|c|c|c|}
\hline
$\times$ &t&x&y&z \\
\hline
t & $t^{2}$ & zt & xt & yt \\
\hline
x&xt&$x^{2}$&xy&zt\\
\hline
y&yt&xt&$y^{2}$&xy\\
\hline
z &zt&xy&yt&$z^{2}$\\
\hline
\end{tabular}
  \hspace{5cm}
\begin{tabular}{|c|c|c|c|c|}
\hline
$\times$ &t&x&y&z \\
\hline
$t$ & $t^{2}$ & zt & xt & yt \\
\hline
$x$&xt&yt&$t^{2}$&zt\\
\hline
$y$&yt&xt&zt&$t^{2}$\\
\hline
$z$ &zt&$t^{2}$& yt&xt\\
\hline
\end{tabular}\\
\qquad\mbox{Table 2: any square is different from the products in}
 (\ref{produits1})
\qquad \qquad \quad  \mbox{Table 3}\\

The case of Table 2  corresponds for instance to the group
$SL(2,\Z /3\Z)$ of the $2\times 2$ matrices with ceofficients in
$\Z /3\Z$, with any of its four-elements conjugacy classes, while
the case of Table 3 corresponds for instance to the alternating
group ${\cal A}_{4}$ of order 12, with any of its four-elements
conjugacy classes. To explicitly solve these nonlinear systems
(\ref{regularcase2}) and (\ref{regularcase1}) one needs more
precision on the group $G$. We solve system (\ref{regularcase2})
in detail in Section 4 for $A_4$. There is in fact a fundamental
difference between the two case, for instance the connection
corresponding in (\ref{module1}) to $\alpha=\beta =\gamma=\lambda$
is a solution of (\ref{regularcase2}) but not a solution of
(\ref{regularcase1}).

We also want to find the `Levi-Civita connection', namely a
regular torsion free and cotorsion free connection for a natural
metric. We need for that end to find a suitable coframing or
metric. As shown in \cite{maj} a natural choice in the group or
quantum group case is to take any Ad-invariant nondegenerate
bilinear form $\eta$ defined on $\Omega_{0}^*$, and indeed
\cite{maj} provides a general `braided-Killing form' construction
that can achieve this. Our `cyclic' conjugacy class ${\cal C}$
described above is not, however, semi-simple in the sense of
Prop.5.4 of \cite{maj} (the braided-Killing form is degenerate)
and we instead have to determine all possible $\eta$.

\begin{theo}\label{theorem3.1}
Let  ${\cal C}$ be a cyclic conjugacy class with $n$ elements.
Then up to normalisation, all nondegenerate $\Ad$-invariant
bilinear forms on $\Omega_0^*$ are given by
\begin{eqnarray}\label{eta}
 \eta^{a,b}=\delta_{a,b}+\mu
\end{eqnarray}
for a
constant $\mu\ne -1/n$. The associated metric in  the
Maurer-Cartan framing is
\[g=\sum_{a\in \CC} e_a\tens_H e_a+\mu
\theta\tens_H\theta.\]
\end{theo}
\proof:

 Here $g$ corresponds to an element $\eta
\in\Omega_0\tens\Omega_0$ with coefficients $\eta^{a,b}$. We
require it to be $\Ad$-invariant and for the matrix of
coefficients to be invertible (this is said more abstractly in
\cite{maj} to handle the quantum group case). The first condition
is easily seen to be the requirement
\begin{eqnarray}\label{adinvariance}
\eta ^{g^{-1}ag,b}=\eta ^{a,gbg^{-1}},\quad \forall a,b\in\CC,
\quad  g\in G.
\end{eqnarray}
This and nondegeneracy is easy to see for the $\eta$ as stated.

Conversely, let us suppose that $\eta$ is Ad-invariant and show
that it is of the form (\ref{eta}). Since $\eta$ is Ad-invariant,
it obeys (\ref{adinvariance}).  We assume the existence of $t\in
\CC$ as in Definition ~\ref{cyclicclass}, then $\Ad_{t}$ is a
cyclic permutation of ${\cal C}-\{t\}$. {}From invariance (\ref
{adinvariance}) it is obvious that 
$\eta^{t,a}=\eta^{t,\Ad_t(a)}=\eta^{t,\Ad_{t^2}(a)}=\cdots
=\eta^{t,\Ad_{t^{n-2}}(a)}$ for any $a\ne t$, and hence by
cyclicity 
\[ \eta^{t,b}=\mu_1,\quad \forall b\ne t,\] for some
constant $\mu_1$. But also from cyclicity we know that for any
$a\in \CC$ there is an element $c\in\CC$ such that $a=ctc^{-1}$.
Hence from (\ref{adinvariance}) we also have
\begin{eqnarray}\label{eachlineeta}
\eta^{a,b}=\eta^{ctc^{-1},b}=\eta^{t,c^{-1}bc}=\mu_1,\quad \forall
a\ne b,
\end{eqnarray}
so all off-diagonals are $\mu_1$. Similarly, we have
$\eta^{a,a}=\eta^{ctc^{-1},ctc^{-1}}=\eta^{t,t}=\mu_2$ for all
$a\in \CC$ by $\Ad$-invariance, for some constant $\mu_2$. Thus
$\eta^{a,b}=(\mu_{2}-\mu_{1})\delta_{a,b}+\mu_{1}$, which has, up
to an overall scaling,  the form stated. The remaining condition
on the parameter $\mu$ comes from the fact that $\eta$ is
invertible. Finally, given $\eta$ we define \[
e^{*a}=\sum_{b\in\CC}e_b\eta^{ba}\] as explained in \cite{maj} for
the associated coframing, which corresponds to the metric $g$ as
stated. \eproof

One can then observe that this  metric is symmetric in the sense
\begin{eqnarray}
\wedge g=0
\end{eqnarray}
The groups $S_{3}$, $SL(2,\Z/3\Z)$ and ${\cal A}_{4}$ are the
examples of groups which obey the hypothesis of the previous
theorem. The theorem clarifies the observation in \cite{maj} for
$S_3$ where $\eta^{a,b}=\delta^{a,b}$ is derived as the braided
Killing form (up to a normalisation) but it is explained that one
may add a multiple $\mu \theta\tens_H\theta$ to the metric
(without changing the connection and Riemmanian curvature). We are
now ready to describe torsion free and cotorsion free connections
in our cyclic case.

\begin{prop}\label{prop3.4}
In the setting of Propositions~\ref{prop3.1} and~\ref{prop3.2} and
for the coframing given by $\eta$ as above, the  torsion free and
cotorsion free connections obey the following relations:
\begin{eqnarray}\label{torandcotor}
R_{t}^{-1}(\alpha)&=&R_{x}^{-1}(\lambda)=R_{y}^{-1}(\beta)
=R_{z}^{-1}(\gamma)\\\nonumber
R_{t}^{-1}(\lambda)&=&R_{x}^{-1}(\alpha)=R_{y}^{-1}(\gamma)
=R_{z}^{-1}(\beta)\\ \nonumber
R_{t}^{-1}(\beta)&=&R_{x}^{-1}(\gamma)=R_{y}^{-1}(\alpha)
=R_{z}^{-1}(\lambda)\\ \nonumber
R_{t}^{-1}(\gamma)&=&R_{x}^{-1}(\beta)=R_{y}^{-1}(\lambda)
=R_{z}^{-1}(\alpha)\\ \nonumber
\end{eqnarray}
where $\alpha,\beta,\gamma,\lambda $ are as in Proposition
\ref{prop3.2}.
\end{prop}
\proof:\\
As in \cite{maj}, when the coframing is given by the
framing and an $\Ad$-invariant $\eta$, one may easily compute the
form of the cotorsion. One has,
\[
D_{A}e^{*a}=\eta^{ba}\extd e_{b}+\sum\limits_{b\in {\cal C},\;c\in
{\cal C}} \eta^{ba}e_{cbc^{-1}}\wedge A_{c}-
 \sum\limits_{b, c \in {\cal C}}\eta^{ba}e_{b}\wedge
 A_{c}\]
as a special case of the quantum groups computation in \cite{maj}.
Since we suppose the connections to be torsion free, equation
(\ref{sommedesAi})  holds, then (cancelling $\eta^{ba}$),
vanishing of cotorsion in equation (\ref{cotorsion}) can be
written equivalently as
\begin{eqnarray}\label{equacortor2}
\extd e_{a}+\sum\limits_{b \in {\cal C}}e_{bab^{-1}}\wedge
 A_{b}=0
\end{eqnarray}
$\forall a \in {\cal C}$. If we write  equation
(\ref{equacortor2}) for $a=t,x,y,z$ respectively, using Table~1,
equations (\ref{les2-form}) and the definition of $\eta $, we
obtain the following system of equations:
 \begin{eqnarray}\label{systemtoretcotor}
e_{t}\wedge A_{t}&+&e_{x}\wedge A_{z}+ e_{y}\wedge
A_{x}+e_{z}\wedge A_{y}-e_{x}\wedge e_{z}-e_{y}\wedge
e_{x}-e_{z}\wedge e_{y}=0\\ \nonumber e_{t}\wedge
A_{z}&+&e_{x}\wedge A_{t}+ e_{y}\wedge A_{y}+e_{z}\wedge
A_{x}-e_{z}\wedge e_{x}-e_{x}\wedge e_{t}-e_{t}\wedge
e_{z}=0\\\nonumber e_{t}\wedge A_{x}&+&e_{x}\wedge A_{y}+
e_{y}\wedge A_{t}+e_{z}\wedge A_{z}-e_{t}\wedge e_{x}-e_{x}\wedge
e_{y}-e_{y}\wedge e_{t}=0\\\nonumber e_{t}\wedge
A_{y}&+&e_{x}\wedge A_{x}+ e_{y}\wedge A_{z}+e_{z}\wedge
A_{t}-e_{z}\wedge e_{t}-e_{t}\wedge e_{y}-e_{y}\wedge e_{z}=0
\end{eqnarray}
we can get the fourth equation of system (\ref{systemtoretcotor})
from the three other. We then  solve only the first three
equations of this system. For that end, we set
\begin{eqnarray}\label{conecadroite}
A_{a}=e_{b}A_{a}^{\prime b}
\end{eqnarray}
 $\forall a \in {\cal C}$, with summation understood for $b\in {\cal
   C}$ ,
and where we set
\begin{eqnarray}\label{diagdeA}
A_{t}^{\prime t}=1+\alpha^{\prime}\quad
 A_{x}^{\prime x}=1+\beta^{\prime}\quad
A_{y}^{\prime y}=1+\gamma^{\prime}\quad A_{z}^{\prime
z}=1+\lambda^{\prime}\quad
\end{eqnarray}
as above. We then proceed in the same manner as we solved sytem
(\ref{syseqtorsion}), using this time the right module structure
of $\Omega^{1}(H)$. We find that the solutions $(A_{a})$ take the
form
\begin{eqnarray}\label{solutiontorcotor}
A_{t}&=&e_{t}(1+\alpha^{\prime})+e_{x}\lambda^{\prime}
+e_{y}\beta^{\prime}
+e_{z}\gamma^{\prime}\\ \nonumber
A_{x}&=&e_{t}\gamma^{\prime}+e_{x}(1+\beta^{\prime})
+e_{y}\lambda^{\prime}
+e_{z}\alpha^{\prime}\\ \nonumber
A_{y}&=&e_{t}\lambda^{\prime}+e_{x}\alpha^{\prime}+
e_{y}(1+\gamma^{\prime}) +e_{z}\beta^{\prime}\\ \nonumber
A_{z}&=&e_{t}\beta^{\prime}+e_{x}\gamma^{\prime}+e_{y}\alpha^{\prime}
+e_{z}(1+\lambda^{\prime})\\ \nonumber &\mbox{with}&\quad
\alpha^{\prime}+\beta^{\prime}+\gamma^{\prime}+\lambda^{\prime}=-1,
\quad
\alpha^{\prime},\beta^{\prime},\gamma^{\prime},\lambda^{\prime}\in
H.
\end{eqnarray}
we then write these solutions using the left module structure on
1-form via the commutation relation in (\ref{extdif}). And we
compare the result to that of system (\ref{module1}) to obtain
system (\ref{torandcotor}) as stated.
\eproof \\

At this level, we get many torsion free and cotorsion free
connections. As one can remark, these equations for the connection
do not depend on the coefficient $\mu$ of $\theta\tens_H\theta$,
just as was the case for $S_3$ in \cite{maj}. Modulo these modes,
we see that there is an essentially unique form of invariant
metric on $G$ and we have given some conditions for the associated
moduli of torsion free and cotorsion free regular connections.

\section{Riemannian geometry of ${\cal A}_{4}$  }

In this section we specialise to the group ${\cal A}_{4}$ and
present stronger results that depend on its structure and not only
on the cyclic form of the conjugacy class. The group is defined by
\begin{eqnarray}\label{groupalterne}
{\cal A}_{4}=\{e,u,v,w,t,x,y,z,t^{2},ut^{2},vt^{2},wt^{2}\}
\end{eqnarray}
where $e$ is the group identity (this should not be confused with
the Maurer-Cartan 1-form) and $t,u,v,w$ are the following
permutations of $\{1,2,3,4\}$:
\begin{eqnarray}\label{generateurA4}
t=(123),\quad u=(14)(23),\quad v=(12)(34),\quad w=(13)(24)
\end{eqnarray}
and
\begin{eqnarray}\label{xyz}
x=tv=ut=(134),\quad y=tw=vt=(243),\quad z=tu=wt=(142).
\end{eqnarray}
The other products are
\begin{eqnarray}\label{produitdansA4}
v^{2}=e,\;w^{2}=e,\;u^{2}=e,\;t^{3}=e,\;
vw=wv=u,\;vu=uv=w,\;wu=uw=v
\end{eqnarray}
and we choose the conjugacy class
\begin{eqnarray}\label{classachoisie}
{\cal C}=\{t,x,y,z\},
\end{eqnarray}
which is `cyclic'. Indeed, we have
\begin{eqnarray*}
\Ad_{t}(x)=txt^{2}=t(ut)t^{2}=z,\;
\Ad_{t}(y)=tyt^{2}=t(vt)t^{2}=x\; \Ad_{t}(z)=t(wt)t^{2}=y
\end{eqnarray*}
and
\begin{eqnarray*}
\Ad_{x}(t)=xtyt=y,\;
Ad_{y}(t)=ytzt=z,\;\Ad_{z}(t)=ztxt=x,\;Ad_{t}(t)=t
\end{eqnarray*}
which show that ${\cal C}$  obeys conditions of
Definition~\ref{cyclicclass}.

One may  also check that the multiplication table of ${\cal C}$
corresponds to Table~3 as announced, so that we have at least one
regular torsion free and cotorsion free connection on the bundle
$H\otimes H$ where $H$ denotes from now $\C [{\cal A}_{4}]$.

\begin{prop}\label{prop4}
For the cyclic conjugacy class on $A_4$, framing defined by the
Maurer-Cartan form and coframing $e^{*}$ by $\Ad$-invariant
$\eta$, there exists a unique `Levi-Civita' connection with
component 1-forms
\begin{eqnarray}\label{levicita}
A_{a}=e_{a}-\frac {1}{4}\theta,\quad \forall a\in {\cal C}
\end{eqnarray}
\end{prop}
\proof: \\
The connection defined in (\ref{levicita}) is easily seen to be a
solution of systems (\ref{module1}) and (\ref{solutiontorcotor}).
We are going to show that it is the unique
 torsion free and cotorsion free connection which is solution of
system (\ref{regularcase2}). Using the properties of operators
$(R_{g})_{g\in {\cal A}_{4}}$ and equations of system
(\ref{torandcotor} ) we find that
\begin{eqnarray}\label{relationentrepara}
\alpha=R_{u}(\lambda),\quad \beta=R_{w}(\lambda),\quad
\gamma=R_{v}(\lambda),
\end{eqnarray}
where $\lambda$ is any function on ${\cal A}_{4}$ which obeys
\begin{eqnarray}\label{condisurlambda}
(R_{u}+R_{v}+R_{w}+\id)(\lambda)=-1.
\end{eqnarray}
At this level, $\lambda$ is not necessarily a scalar. To determine
it, we set
\begin{eqnarray}\label{lamdaencordonnee}
\lambda=\sum\limits_{{g\in {\cal A}_{4}}}\lambda_{g}\delta_{g}
\end{eqnarray}
hence we get from (\ref{condisurlambda}) that
\begin{eqnarray*}\label{lamdaredui}
\lambda&=&(-1-\lambda_{v}-\lambda_{w}-\lambda_{u})\delta_{e}+
\lambda_{v}\delta_{v}+
\lambda_{w}\delta_{w}+\lambda_{u}\delta_{u}+(-1-\lambda_{x}
-\lambda_{y}-\lambda_{z}) \delta_{t}
+\lambda_{x}\delta_{x}\\\nonumber
&+&\lambda_{y}\delta_{y}+\lambda_{z}\delta_{z} +
(-1-\lambda_{tx}-\lambda_{ty}-\lambda_{tz})\delta_{t^{2}}+
\lambda_{tx}\delta_{tx}+\lambda_{ty}\delta_{ty}+
\lambda_{tz}\delta_{tz}.
\end{eqnarray*}
We write out the first equation of system (\ref{regularcase2}) in
the basis (\ref{base2-form})  of $\Omega^{2}(H)$, passing from the
right module structure to the left one, then set to zero each
coefficient of the basis element and  obtain the following
equations
\begin{eqnarray}\label{sysfinalderegu}
(1+\alpha)R_{1}(\gamma)&+&\lambda R_{1}(\lambda)+\beta
R_{1}(\alpha)+ \gamma R_{1}(1+\beta)\\ \nonumber &-&\beta
R_{4}(1+\alpha)-\gamma R_{4}(\beta)-\alpha
R_{4}(\gamma)-(1+\lambda)R_{4}(\lambda)=0 \\ \nonumber
(1+\alpha)R_{1}(\lambda)&+&\lambda R_{1}(1+\gamma)+\beta
R_{1}(\beta)+ \gamma R_{1}(\alpha)\\ \nonumber &-&\lambda
R_{3}(\gamma)-\alpha R_{3}(\lambda)-(1+\gamma) R_{3}(\alpha)-\beta
R_{3}(1+\beta)=0 \\ \nonumber (1+\alpha)R_{1}(\beta)&+&\lambda
R_{1}(\alpha)+\beta R_{1}(1+\lambda)+ \gamma R_{1}(\gamma)\\
\nonumber &-&\beta R_{4}(\lambda)-\gamma R_{4}(1+\gamma)-\alpha
R_{4}(\beta)-(1+\lambda)R_{4}(\alpha)=0 \\ \nonumber \gamma
R_{2}(1+\alpha)&+&(1+\beta) R_{2}(\beta)+\lambda R_{2}(\gamma)+
\alpha R_{2}(\lambda)\\ \nonumber &-&\lambda R_{3}(\gamma)-\alpha
R_{3}(\lambda)-(1+\gamma) R_{3}(\alpha)-\beta R_{3}(1+\beta)=0 \\
\nonumber \gamma R_{2}(\lambda)&+&(1+\beta)
R_{2}(1+\gamma)+\lambda R_{2}(\beta)+ \alpha R_{2}(\alpha)\\
\nonumber &-&\beta R_{4}(\gamma)-\gamma R_{4}(\lambda)-\alpha
R_{4}(\alpha)-(1+\lambda) R_{4}(1+\beta)=0 \\ \nonumber \gamma
R_{2}(\beta)&+&(1+\beta) R_{2}(\alpha)+\lambda R_{2}(1+\lambda)+
\alpha R_{2}(\gamma)\\ \nonumber &-&\beta R_{4}(1+\alpha)-\gamma
R_{4}(\beta)-\alpha R_{4}(\gamma)-(1+\lambda) R_{4}(\lambda)=0 \\
\nonumber \lambda R_{3}(1+\alpha)&+&\alpha R_{3}(\beta)+(1+\gamma)
R_{3}(\gamma)+ \beta R_{3}(\lambda)\\ \nonumber &-&\beta
R_{4}(\lambda)-\gamma R_{4}(1+\gamma)-\alpha
R_{4}(\beta)-(1+\lambda) R_{4}(\alpha)=0 \\ \nonumber \lambda
R_{3}(\beta)&+&\alpha R_{3}(\alpha)+(1+\gamma) R_{3}(1+\lambda)+
\beta R_{3}(\gamma)\\ \nonumber &-&\beta R_{4}(\gamma)-\gamma
R_{4}(\lambda)-\alpha R_{4}(\alpha)-(1+\lambda) R_{4}(1+\beta)=0,
\end{eqnarray}
where the indexes $1,2,3,4$ refer respectively to $t,x,y,z$. We
then use (\ref{relationentrepara}) to write these equations
respectively in terms of $\lambda$,  then in terms of its scalar
components. A long but straightforward  computation of these
components leads to $\lambda_{g}=-\frac{1}{4}, \forall g \in {\cal
A}_{4} $, hence as an element of $H$, $\lambda=-\frac{1}{4}$. {}From
(\ref{relationentrepara}), we also have
$\alpha=\beta=\gamma=-\frac{1}{4}$. The expression of the
corresponding connection in (\ref{module1}) is then as stated. To
end the proof of the Proposition~\ref{prop4}, one checks easily
that this connection is also solution of the other equations of
system (\ref{regularcase2}).
\eproof \\

We refer to this connection as {\em the} `Levi-Civita connection'
for the invariant metric on the group ${\cal A}_{4}$.

\begin{prop}\label{prop5} \quad The covariant derivative $\nabla
:\Omega^{1}(H) \longrightarrow
 \Omega^{1}(H)\mH \Omega^{1}(H)$ for
 the above `Levi-Civita connection'
on ${\cal A}_{4}$, and its Riemann curvature  ${\cal R}
:\Omega^{1}(H) \longrightarrow
 \Omega^{2}(H)\mH \Omega^{1}(H)$ are given by
\begin{eqnarray}\label{covarianderivative}
\nabla(e_{t})&=&-e_{t}\mH e_{t}-e_{x}\mH e_{z}-e_{y}\mH
e_{x}-e_{z}\mH e_{y}+{1\over 4}\theta \mH \theta \\ \nonumber
\nabla (e_{x})&=&-e_{t}\mH e_{y}-e_{x}\mH e_{x}-e_{y}\mH
e_{z}-e_{z}\mH e_{t} +{1\over 4}\theta \mH \theta\\ \nonumber
\nabla(e_{y})&=&-e_{t}\mH e_{z}-e_{x}\mH e_{t}-e_{y}\mH
e_{y}-e_{z}\mH e_{x} +{1\over 4}\theta \mH \theta\\ \nonumber
\nabla (e_{z})&=& -e_{t}\mH e_{x}-e_{x}\mH
e_{y}-e_{y}\mH e_{t}-e_{z}\mH e_{z}+{1\over 4}\theta \mH \theta.\\
\nonumber
\end{eqnarray}
\begin{eqnarray}\label{formedecourbure}
{\cal R}(e_{t})&=&\extd e_{t}\mH e_{t}+\extd e_{x}\mH e_{z}+\extd
e_{y}\mH e_{x}+\extd e_{z}\mH e_{y}\\ \nonumber {\cal
R}(e_{x})&=&\extd e_{t}\mH e_{y}+\extd e_{x}\mH e_{x}+\extd
e_{y}\mH e_{z}+\extd e_{z}\mH e_{t}\\ \nonumber {\cal
R}(e_{y})&=&\extd e_{t}\mH e_{z}+\extd e_{x}\mH e_{t}+\extd
e_{y}\mH e_{y}+\extd e_{z}\mH e_{x}\\ \nonumber {\cal R}(e_{z})&=&
\extd e_{t}\mH e_{x}+\extd e_{x}\mH e_{y}+\extd e_{y}\mH e_{t}
+\extd e_{z}\mH e_{z}\\
\nonumber
\end{eqnarray}
\end{prop}
\proof: The curvature 2-form $F$ is defined  by equation
(\ref{2formdecourbure}). In the present case, we have $bc\notin
{\cal C},\forall b,c\in {\cal
  C}$,
so that $\sum\limits_{b,c\in {\cal C},bc=a}A_{b}\wedge
A_{c}=0$,$\forall a \in {\cal C} $. We have also
$\sum\limits_{a\in {\cal C}}A_{a}=0$ and $\extd\theta=0$, hence
$F_{a}=\extd A_{a}=\extd e_a$ for the form of the connection in
(\ref{levicita}). This is exactly the same argument as for $S_3$
in \cite{maj}. Next, if we replace $\alpha$ in  formula
(\ref{formuledelacourbureriem}) by $e_{t},e_{x},e_{y},e_{z}$
respectively, and use Table~1, we obtain relations
(\ref{formedecourbure}) for the curvature. Finally, we compute the
value of the covariant derivative on the basis 1-forms $\{e_{a}\}$
by using  formula (\ref{covariantderivative}). Explicitly, we have
\begin{eqnarray*}
\nabla (e_{a})&=&-\sum\limits_{b\in{\cal
    C}}A_{b}\mH(e_{b^{-1}ab}-e_{a})
\\ \nonumber
&=&-\sum\limits_{b\in{\cal C}}(e_{b}-{1\over 4}\theta)\mH
(e_{b^{-1}ab}-e_{a})
\\ \nonumber
&=&-\sum\limits_{b\in{\cal C}}e_{b}\mH (e_{b^{-1}ab}-e_{a})+
{1\over 4}\theta \mH \sum\limits_{b\in{\cal C}}
(e_{b^{-1}ab}-e_{a})
 \\ \nonumber
&=&-\sum\limits_{b\in{\cal C}}e_{b}\mH e_{b^{-1}ab}+
\sum\limits_{b\in{\cal C}}  e_{b}\mH e_{a}+{1\over 4} \theta \mH
\sum\limits_{b\in{\cal C}} (e_{b}-e_{a})
\\ \nonumber
&=&-\sum\limits_{b\in{\cal C}}e_{b}\mH e_{b^{-1}ab}+{1\over 4}
\theta\mH \theta
\end{eqnarray*}
According to Table~1, this last equation gives relations
(\ref{covarianderivative}) as stated .
\eproof \\

{}From the Riemann curvature and the canonical lift $i$ we can
compute the Ricci curvature of the Levi-Civita connection on $A_4$
and find that it vanishes. In fact we can prove a slightly
stronger result that is it the {\em only} Ricci flat connection
for this choice of framing.

\begin{theo} For the framing defined by the
Maurer-Cartan 1-form, and for the canonical lift $i$, the above
Levi-Civita connection on $A_4$ is the unique regular Ricci-flat
connection.
\end{theo}
\proof:\\
In the present case, the canonical lift takes the form
\begin{eqnarray}\label{produitinterieur2}
i(e_{a}\wedge e_{b})=e_{a}\mH e_{b}-
\frac{1}{3}\sum\limits_{cd=ab,c\neq d}e_{c}\mH e_{d}, \quad
i(e_{a}\wedge e_{a})=0.
\end{eqnarray}
We have to solve for vanishing of\cite{maj}
\[Ricci=\sum\limits_{a\in {\cal
C}}\<f^{a}, (i\mH \id){\cal
  R}(e_{a})\>=\sum\limits_{a,b,c \in {\cal C}}i(F_c)^{ab}e_{b} \mH
(e_{c^{-1}a c}-e_{a}),\] where $i(F_{c})=i(F_c)^{ab}e_{a}\mH
e_{b}$, and the pairing is made between each $f^{a}$ and the first
factor of the tensor product $(i\mH \id){\cal
  R}(e_{a})$ according to the formula
$\<f^{a},me_{b}\>=m\delta_{b}^{a},\forall m \in H$. In our case
this becomes
\begin{eqnarray}\label{ricciA4vanishing}
&{}&\<f^{t},i(F_{x})\mH (e_{z}-e_{t})+i(F_{y})\mH (e_{x}-e_{t})+
i(F_{z})\mH (e_{y}-e_{t})\>\\ \nonumber &+& \<f^{x},i(F_{t})\mH
(e_{y}-e_{x})+i(F_{y})\mH (e_{z}-e_{x})+ i(F_{z})\mH
(e_{t}-e_{x})\>\\ \nonumber &+& \<f^{y},i(F_{t})\mH
(e_{z}-e_{y})+i(F_{x})\mH (e_{t}-e_{y})+ i(F_{z})\mH
(e_{x}-e_{y})\>\\ \nonumber &+& \<f^{z},i(F_{t})\mH
(e_{x}-e_{z})+i(F_{x})\mH (e_{y}-e_{z})+ i(F_{y})\mH
(e_{t}-e_{z})\>=0
\end{eqnarray}
We first compute $F_{t},\;F_{x},\;F_{y},\;F_{z}$ and
$i(F_{t}),\;i(F_{x}),\; i(F_{y}),\; i(F_{z})$ for general free
torsion connections given in Proposition~\ref{prop3.2}, then we
rewrite equation (\ref{ricciA4vanishing})
 in terms of  the basic elements
$\{e_{a}\mH e_{b}\}_{a,b \in {\cal C}}$ of the left H-module
$\Omega^{1}(H) \mH \Omega^{1}(H)$, the vanishing of each
coefficient of the mentioned basic elements leads to 16 equations
in terms of $\alpha,\beta,\gamma, \lambda$ and their `first order
derivatives'
$\partial^{a}\alpha,\partial^{a}\beta,\partial^{a}\gamma,\partial^{a}
\lambda$,\quad $a\in {\cal C}$.

We find that it is enough to solve the following 4 equations
coming from the coefficients of $e_{t}\mH e_{t},\; e_{x}\mH
e_{x},\; e_{y}\mH e_{y}, \;e_{z}\mH e_{z}$ respectively:
\begin{eqnarray}
\beta&-&2\gamma+\lambda+\partial^{y}\alpha+\partial^{z}\alpha
+\partial^{x}\alpha
+\partial^{t}\beta-2\partial^{z}\beta+\partial^{t}\gamma
-2\partial^{x}\gamma
+\partial^{t}\lambda-2\partial^{y}\lambda=0 \\ \nonumber
\alpha&-&3\beta+\lambda+\gamma+\partial^{x}\alpha
-2\partial^{y}\alpha
+\partial^{y}\beta+\partial^{t}\beta+\partial^{z}\beta+
\partial^{t}\gamma -2\partial^{z}\gamma+
\partial^{x}\lambda-2\partial^{t}\lambda=0 \\ \nonumber
\alpha
&+&\beta-3\gamma+\lambda+\partial^{y}\alpha-2\partial^{z}\alpha
+\partial^{y}\beta-2\partial^{t}\beta+
\partial^{x}\gamma+\partial^{t}\gamma+\partial^{z}\gamma
+\partial^{y}\lambda-2\partial^{x}\lambda=0\\ \nonumber \alpha
&+&\beta+\gamma-3\lambda +\partial^{z}\alpha -2\partial^{x}\alpha
+\partial^{z}\beta -2\partial^{y}\beta +\partial^{z}\gamma
-2\partial^{t}\gamma
+\partial^{y}\lambda+\partial^{t}\lambda+\partial^{x}\lambda=0
\end{eqnarray}
Indeed, we transform these 4 equations to a system of 48 linear
equations where the variables are the components of
$\alpha,\beta,\gamma$ and $\lambda$ in the basis
$(\delta_{g})_{g\in {\cal A}_{4}}$. The unique solution of the
mentioned system which obeys the condition
$\alpha+\beta+\gamma+\lambda=-1$ as in Proposition \ref{prop3.2}
is $\alpha=\beta=\gamma=\lambda=-1/4$. To end the proof one just
checks easily that this solution is also a solution of the 12
remaining equations (of the 16 ones mentioned above), coming from
the coefficients of $e_{a}\mH e_{b},\; a\neq b$ in equation
(\ref{ricciA4vanishing}). \eproof

One can also check that the Ricci tensor for the Levi-Civita
connection with respect to the alternative `lift'
\begin{eqnarray}\label{produitinterieur1}
i'(e_{a}\wedge e_{b})=e_{a}\mH e_{b}-e_{aba^{-1}}\mH e_{a},
\end{eqnarray}
also vanishes, i.e. the result does not depend strongly on the
choice of lift. This is the same as found for $S_{3}$, where the
two Ricci tensors with respect to $i$ and $i'$ respectively are
the same up to a scale \cite{maj}.

\section{The Dirac operator for ${\cal A}_{4}$}

Following the formalism of reference \cite{maj}, we write down in
this section the `gamma matrices' and the Dirac operator
associated to the Maurer-Cartan framing $e$ and the coframing
$e^{*}$ for the invariant metric. We use the associated
Levi-Civita connection constructed above.

 For the `spinor'
representation, we consider the standard 3-dimensional
representation of ${\cal A}_{4}$ defined on a vector space $W$ by
\begin{eqnarray}\label{representation}
\rho_{W}(t)&=&\left (\begin{array}{ccc}
0&0&1\\
1&0&0\\
0&1&0
\end{array}\right),
 \quad
\rho_{W}(u)=\left (\begin{array}{ccc}
-1&0&0\\
0&-1&0\\
0&0&1
 \end{array}\right),\\ \nonumber
\quad \rho_{W}(v)&=&\left (\begin{array}{ccc}
-1&0&0\\
0&1&0\\
0&0&-1
\end{array}\right),
\quad \rho_{W}(w)=\left (\begin{array}{ccc}
1&0&0\\
0&-1&0\\
0&0&-1
\end{array}\right)
\end{eqnarray}
and $e$ is the unit matrix $I$. The Casimir element $C$ associated
to the operator $\eta$ is given in \cite{maj} by
\begin{eqnarray}
C=\eta_{ab}^{-1}f^{a}f^{b}=\eta_{ab}^{-1}(a-e)(b-e)
\end{eqnarray}
with summation understood, $b,a \in {\cal C}$. One checks that it
corresponds in the general case of the class ${\cal
C}=\{t,x,y,z\}$ as in Prop.\ref{prop3.1}, to the explicit form
\begin{eqnarray}\label{casimir}
C&=&\frac{1+3\mu}{1+4\mu}[t^{2}+x^{2}+y^{2}+z^{2}-2(t+x+y+z)+4e]\\ 
\nonumber
&+&\frac{-3\mu}{1+4\mu}[tx+ty+tz+xy-2(t+x+y+z)+4e]
\end{eqnarray}

 In the case of ${\cal A}_{4}$,
equation (\ref {casimir}) reads
\begin{eqnarray*}\label{casimirdeA4}
C=\frac{1}{1+4\mu}[(t-e)^{2}+(x-e)^{2}+(y-e)^{2}+(z-e)^{2}]
\end{eqnarray*}
then
\begin{eqnarray*}
\rho _{W}(C)=\frac{4}{1+4\mu}I.
\end{eqnarray*}
Next, we choose our gamma-matrix $\gamma$  to be the `tautological
gamma-matrix' \cite{maj}  associated to $\rho_{W}$ and $\eta$
defined by
\begin{eqnarray}
\gamma_{a}=\eta_{a b}^{-1}\rho _{W}(f^{b})=\sum\limits_{b \in{\cal
C}}\eta_{a b}^{-1}\rho_{W}(b-e),\quad \forall a\in\CC.
\end{eqnarray}
In our case we find
\begin{eqnarray}\label{simpledefdegamma}
\gamma_{a}=\rho _{W}(a-e)+\frac{4\mu}{1+4\mu}
\end{eqnarray}
and that these matrices  obey the relations
\begin{eqnarray}\label{sommedesgamma}
\sum\limits_{a \in {\cal C}}\gamma_{a}=-\frac{4}{1+4\mu},
\end{eqnarray}
\begin{eqnarray}\label{commutaiondesgamma}
\gamma_{a}\gamma_{b}+\gamma_{b}\gamma_{a}+
\frac{2}{1+4\mu}(\gamma_{a}+\gamma_{b})+\frac{2}{(1+4\mu)^{2}} 
=\rho_{W}(a b
+b a)
\end{eqnarray}
following directly from (\ref{simpledefdegamma}).

Equations (\ref{simpledefdegamma}), (\ref{sommedesgamma}) and
(\ref{commutaiondesgamma}) hold in the general case considered in
Proposition~\ref{prop3.1}, providing that the multiplication's
table is that of Table~3. The explicit matrix representation of
these gamma-matrices above for ${\cal A}_{4}$ are:

\begin{eqnarray}\label{repreexplicitedesgamma}
\gamma_{t}&=&\left (\begin{array}{ccc}
\frac{-1}{1+4\mu}&0&1\\
1&\frac{-1}{1+4\mu}& 0\\
0&1&\frac{-1}{1+4\mu}
\end{array}\right),\qquad 
\gamma_{x}=\left (\begin{array}{ccc}
\frac{-1}{1+4\mu}&0&-1\\
-1&-\frac {-1}{1+4\mu}&0\\
0&1&\frac{-1}{1+4\mu}
 \end{array}\right),\\ \nonumber
\quad \gamma_{y}&=&\left (\begin{array}{ccc}
\frac{-1}{1+4\mu}&0&-1\\
1&\frac {-1}{1+4\mu}&0\\
0&-1&\frac{-1}{1+4\mu}
\end{array}\right),\qquad
\gamma_{z}=\left (\begin{array}{ccc}
\frac{-1}{1+4\mu}&0&1\\
-1&\frac{-1}{1+4\mu}&0\\
0&-1&\frac{-1}{1+4\mu}
\quad \
\end{array}\right)
\end{eqnarray}

\begin{prop}\label{prop6}
The Dirac operator (\ref{operateurDirac}) on ${\cal A}_{4}$ for the
gamma-matrices and the Levi-Civita connection on ${\cal A}_{4}$
constructed above is given by

\begin{eqnarray}\label{dirac}
\D&=&\partial^{a}\gamma_{a}-4
\end{eqnarray}
 (sum over $ a \in {\cal C}$). For $\mu=0$ we have explicitly
\[ \D=\left (\begin{array}{ccc}
-R_{t}-R_{x}-R_{y}-R_{z}&
0 &R_{t}-R_{x}-R_{y}+R_{z}
\\
 R_{t}-R_{x}+R_{y}-R_{z} & -R_{t}-R_{x}-R_{y}-R_{z} &
 0\\
 0 &R_{t}+R_{x}-R_{y}-R_{z} &
-R_{t}-R_{x}-R_{y}-R_{z}
\end{array}\right).\]
This has 18 zero modes, 3 modes with eigenvalue $\pm 4$, 3
modes with eigenvalue $\pm 4q$, and 3 modes with eigenvalue
 $\pm 4\bar q$, where $q=e^{2\pi\imath/3}$. 
\end{prop}
\proof:\\
The formula giving the Dirac operator in terms of the
gamma-matrices and the  representation $\rho_{W}$  is given by
equation (\ref{operateurDirac}). We  first observe that for the
representation $\rho_{W}$ above, the following two equations hold:
 $\sum\limits_{a \in {\cal C}}\rho_{W}(a)=0$ and
 $\sum\limits_{a \in {\cal C}}\rho_{W}(a ^{2})=0$.
Using the  $A_{a}^{b}$ defined by (\ref{levicita}), and the fact
that
 every element of ${\cal
  C}$ is of order 3, we obtain
\begin{eqnarray*}
\D &=&\partial^{a}\gamma_{a}-\sum\limits_{a,b \in {\cal C}}
(\delta_{a}^{b}-\frac{1}{4})[\rho_{W}(a-e)+\frac{4\mu}{1+4\mu}]\rho_{W}(b
^{2}-e)\\ \nonumber &=&\partial^{a}\gamma_{a}-\sum\limits_{a \in
{\cal C}} [\rho_{W}(a-e)+\frac{4\mu}{1+4\mu}]\rho_{W}(a ^{2}-e)+
\sum\limits_{a \in {\cal
    C}}\frac{1}{4}[\rho_{W}(a-e)+\frac{4\mu}{1+4\mu}](-4I)\\ \nonumber
&=&\partial^{a}\gamma_{a} -\sum\limits_{a \in {\cal C}}
[\rho_{W}(a^{3}-a-a^{2}+e)]-\frac{4\mu}{1+4\mu}\sum\limits_{a \in {\cal
C}}\rho_{W}(a ^{2}-e)-\sum\limits_{a
  \in {\cal C}}\gamma_{a}\\ \nonumber
&=&\partial^{a}\gamma_{a}-4
\end{eqnarray*}
we then replace in equation (\ref{dirac}) the representation of
the gamma-matrices from
 (\ref{repreexplicitedesgamma}) to obtain the matrix representation of $\D$
as stated.

To compute its eigenvalues we need $R_a$ explicitly as 
$12\times 12$ matrices. In the
basis spanned by delta-functions at 
$\{e,u,v,w,t,x,y,z,t^2,ut^2,vt^2,wt^2\}$, the right translation
operators take the form
\[R_t=\left(\begin{array}{ccc}0&I&0\\ 0&0&I\\ I&0&0\end{array}
\right),\quad
R_x=\left(\begin{array}{ccc}0&Y&0\\ 0&0&Z\\ X&0&0\end{array}\right),
\quad
 R_y=\left(\begin{array}{ccc}0&X&0\\ 0&0&Y\\ Z&0&0\end{array}\right),
\quad
R_z=\left(\begin{array}{ccc}0&Z&0\\ 0&0&X\\ Y&0&0\end{array}\right)\]
where $I$ is the $4\times 4$ identity and
\[ X=\left(\begin{array}{cccc}0&0&1&0\\ 0&0&0&1\\ 1&0&0&0\\0&1&0&0
\end{array}\right),\quad
Y=\left(\begin{array}{cccc}0&1&0&0\\ 1&0&0&0\\ 0&0&0&1\\0&0&1&0
\end{array}\right),\quad
Z=\left(\begin{array}{cccc}0&0&0&1\\ 0&0&1&0\\ 0&1&0&0\\1&0&0&0
\end{array}\right).\]
We then obtain the eigenvalues as stated.
\eproof \\

The eigenvalues here do in fact depend on $\mu$ and the case $\mu=0$ seems to
be the more natural since it corresponds 
to the simplest metric $\delta_{a,b}$. The $-4$ in (\ref{dirac}) corresponds 
to the constant curvature
of $A_4$ as for $S_3$ in \cite{maj}. As
for $S_3$, this offset ensures a symmetrical distribution of
eigenvalues about zero.

We will now construct the eigenstates of $\D$. Before doing that 
we look at the spin 0 or scalar wave equation defined by the 
corresponding wave operator 
\eqn{lapl}{ \Box = -\eta^{-1}_{ab}\del^a\del^b=-\sum\limits_{a}
\partial^{a}\partial^{a}=\sum\limits_{a}
(2R_{a}-R_{a^{2}}-\id).}
We do not exactly expect a Lichnerowicz formula relating this to 
the square of $\D$, but we find that it is the square of a first-order 
operator with eigenvalues contained in those of $\D$.
It is  easy to solve the wave equation directly.
\begin{prop}\label{laplacian}
\[ \Box=-{1\over 4}(\sum_a \del^a)^2=-{1\over 4}(D_0-4)^2,\quad 
D_0=\sum_a R_a.\]
There is 1 zero mode, given by the constant function. There is 1 mode 
of eigenvalue $12q$ 
and one of eigenvalue $12\bar q$ given by the two other 1-dimensional 
representations of $A_4$. Finally there are 9 modes with eigenvalue $-4$ 
given by the matrix elements
of the remaining 3-dimensional irreducible representation  $\rho_{W}$ 
above.  
\end{prop}
\proof:\\
The square form of $\Box$ follows from the multiplication Table~3. 
{}From there one finds that $(\sum_a R_a)^2=4\sum_a R_{a^2}$, after 
which the result follows. To solve the wave equation, note that the 
nontrivial 1-dimensional representations $\rho,\bar\rho$ of $A_4$ 
are given by
\[ \rho(t)=q,\quad\rho(u)=\rho(v)=\rho(w)=1,\quad
\bar\rho(t)=\overline{q},\quad\bar\rho(u)=\bar\rho(v)=\bar\rho(w)=1. \] 
Then $\forall m \in {\cal A}_{4}$,
\begin{eqnarray*}
\Box \rho(m) &=& 2\sum\limits_{a}\rho(m)\rho(a)-\sum\limits_{a}
\rho(m)\rho(a^{2})-4\rho(m)\\
&=& ( 8q-4-4\bar{q})\rho(m)=12q\rho(m)
\end{eqnarray*}
Similarly for $\bar\rho$ with $q$ replaced by $\overline{q}$. Finally 
for the matrix elements $\{\rho_{kl}\}$ of $\rho_W$,
we have
\begin{eqnarray*}
\Box \rho_{kl}(m)&=&\sum\limits_{a}
[2R_{a}\rho_{kl}(m)-R_{a^{2}}\rho_{kl}(m)-\rho_{kl}(m)]\\
&=&\sum\limits_{a}\sum\limits_{i}
[2\rho_{ki}(m)\rho_{il}(a)-\rho_{ki}(m)\rho_{il}(a^{2})]
-4\rho_{kl}(m)\\
&=&-4\rho_{kl}(m)
\end{eqnarray*}
since $\sum\limits_{a}\rho_{W}(a^{2})=0$ and
$\sum\limits_{a}\rho_{W}(a)=0$. 
The 9 ``waves'' $\rho_{kl}$ are linearly independent because the
representation
$\rho_{W}$ is irreducible. Hence  we have completely diagonalised 
the 12x12
matrix $\Box$. Equivalently, we have diagonalised $D_0$ with 
corresponding eigenvalues $4,4q,4\bar q,0$ as for $\D$ above.
\eproof

Moreover, every function $\phi$ on ${\cal A}_{4}$ has a
unique decomposition of the form 
\eqn{fou}{\phi=p_{0}+p_{1}\rho+p_{2}\bar\rho
+\sum\limits_{k,l}p_{kl}\rho_{kl}}
for some numbers $p_{0},p_{1},p_{2},p_{kl}$ which are the components
of $\phi$ in the nonabelian Fourier transform. The decomposition above
corresponds precisely to the Peter-Weyl decomposition, just as noted 
for $S_{3}$ in \cite{mr}.

We now use the preceding results to completely solve the Dirac equation. 
We set
\begin{eqnarray*}
 D_{1}=R_{t}-R_{x}-R_{y}+R_{z},\quad
D_{2}&=&R_{t}-R_{x}+R_{y}-R_{z},\quad \; D_{3}=R_{t}+R_{x}-R_{y}-R_{z}
\end{eqnarray*}
so that
\[ \D=\left (\begin{array}{ccc}-D_0&0 &D_1\\D_2 & -D_0 &0\\0 &D_3 &-D_0
\end{array}\right).\]
Let us note first of all that
\begin{eqnarray}
D_{1}D_{2}=0,\;\; D_{2}D_{3}=0,\;\; D_{3}D_{1}=0,\;\; D_{0}D_{i}=0,
\;\; D_{i}^{2}=0,
1\leq i\leq 3,
\end{eqnarray}
from which we see by inspection that the following are 18 
linearly-independent zero modes of $\D$:
\eqn{zerD}{ \left(
\begin{array}{c} D_{2}\rho_{k3}\\ 0 \\ 0 \end{array} \right),
\left(
\begin{array}{c} 0\\ D_{3}\rho_{k1} \\ 0 \end{array} \right),
\left(
\begin{array}{c}  0\\ 0\\ D_{1}\rho_{k2}  \end{array} \right);
\quad \left(
\begin{array}{c}  D_{3}\rho_{k1}\\ 0\\ 0  \end{array} \right),
\left(
\begin{array}{c} 0 \\ D_{1}\rho_{k2}\\ 0 \end{array} \right),
\left(
\begin{array}{c} 0 \\ 0\\ D_{2}\rho_{k3} \end{array} \right)}
for $1\leq k\leq 3$. Similarly it is immediate by inspection that
\eqn{-4qD}{     \left(
\begin{array}{c} \rho^n\\ 0 \\ 0 \end{array} \right),
\left(
\begin{array}{c} 0\\ \rho^n \\ 0 \end{array} \right),
\left(
\begin{array}{c}  0\\ 0\\ \rho^n  \end{array} \right)}
are 3 modes with eigenvalue $-4q^n$, for $n=0,1,2$. This is because 
$D_0\rho=4q \rho$ (as
in Proposition~\ref{laplacian} above) while $D_i\rho=0$ for $i>0$.

It remains only to construct the $+4q^n$ eigenmodes for $n=0,1,2$. 
Before doing this let us make two observations about the modes 
already evident. First of all, let
$\hat\rho$ denote the operator of multiplication by $\rho$. Then 
$R_a\hat\rho=q\hat\rho R_a$ since $\rho(a)=q$ for $a=t,x,y,z$. Hence
\eqn{hatrho}{ \D \hat\rho=q \hat\rho\, \D.}
Thus, multiplication of a spinor mode by the function $\rho$ 
multiplies its eigenvalue by $q$. This
generates the $-4q^n$ modes above from the $n=0$ case. 

Secondly, from the multiplication Table~3 we see that 
\eqn{tD}{ R_t D_1=D_2 R_t}
and its cyclic rotations $1\to 2\to 3\to 1$. {}From this one finds 
\eqn{chi}{ [\chi,\D]=0,\quad \chi^3=\id;\quad\chi
=\left (\begin{array}{ccc}0&0 &R_t\\ R_t & 0 &0\\0 &R_t &0
\end{array}\right).}
This $\chi$ generates the other two modes from the first in 
each group of three in (\ref{zerD}) and (\ref{-4qD}). In the case of 
the zero modes note that 
\eqn{trho}{R_t\rho_{k1}=\rho_{k2}}
and its cyclic rotations $1\to 2\to 3\to 1$, from the explicit 
form of $\rho_W(t)$.

We now observe that if we make an ansatz of the form
\[ \psi=\left(
\begin{array}{c} \phi\\ R_t\phi \\ R_{t^2}\phi \end{array} 
\right)=\left(
\begin{array}{c} \id\\ R_t \\ R_{t^2}\end{array} \right)\phi\]
for function $\phi$ then 
\[ \D\psi=\left(
\begin{array}{c} \id\\ R_t \\ R_{t^2} \end{array} 
\right)(-D_0+R_{t^2}D_2)\phi\]
so eigenspinors are induced by eigenfunctions of the operator 
\[ -D_0+R_{t^2}D_2=-D_0+R_e -R_u-R_v+R_w.\] 
All of the $\rho_{kl}$ are zero modes of $D_0$ (as in 
Proposition~\ref{laplacian}), while among them precisely 
$\rho_{k1}$ is an eigenmode of $R_e-R_u-R_v+R_w$, with eigenvalue 
$4$ (this follows from ${1\over 4}(\rho_W(e)-\rho_W(u)-\rho_W(v)
+\rho_W(w))$ being a projection matrix of rank 1). 
Hence $\phi=\rho_{k1}$ in the ansatz yields three spinor modes
\eqn{4qD}{  \left(
\begin{array}{c} \rho_{k1}\\ \rho_{k2} \\ \rho_{k3} \end{array} 
\right),\quad 1\leq k\leq 3}
with eigenvalue $+4$ of $\D$. Applying $\hat\rho$ generates the 
three with eigenvalue $4q$ and
then the three with eigenvalue $4\bar q$. 

This completes our diagonalisation of $\D$. Finally, we note that 
there necessarily exists an operator $\gamma$ with 
$\gamma^2=\id$ and
$\{\gamma,\D\}=0$, but it is not unique. Thus, diagonalising
$\D$, we can group the eigenbasis into pairs of 3-blocks of zero 
modes according to the two groups in (\ref{zerD}), 
interchanged by $\gamma$, and similarly we define $\gamma$ to 
interchange the two 3-blocks with eigenvalues $\pm 4q^n$. This defines 
at least one choice of $\gamma$, suggested 
by our explicit diagonalisation. 

\section{Cohomology and concluding remarks}

In this paper we have concentrated on the Riemannian geometry of
$A_4$. There are also some more elementary geometrical questions that
one could look at, related to the differential structure alone. We
will discuss some of them here. 

Firstly,  given an exterior algebra $\Omega(H)$ one has a 
noncommutative de Rahm cohomology defined as usual by closed forms 
modulo exact ones. We find, just as
for $S_3$ in \cite{mr}, that $\theta$ generates $H^1$.

\begin{prop} For $A_4$ with cyclic conjugacy class $\{t,x,y,z\}$ the
first noncommutative de Rahm cohomology is
\[ H^1(A_4)=\C.\theta\]
\end{prop}
\proof:\\
We compute $\del^a=R_a-\id$ explicitly as four $12\times 12$ matrices
for their action on $\C[A_4]$. The concatenation of these define 
$\extd_0:\C[A_4]\to \C[A_4]\tens \Omega_0=\Omega^1(H)$ as a 
$48\times 12$-matrix.  We also 
define the $8\times 16$-matrix $\pi$ which sends 
$e_a\tens_He_b\to e_a\wedge e_b$ 
using the tensor product 16-dimensional basis of 
$\Omega_0\tens\Omega_0$ and the 8-dimensional vector space over 
$\C$ with basis (\ref{base2-form}). Similarly, we define an 
$8\times 4$-matrix for $\extd$ acting on $\Omega_0$ again 
using the basis (\ref{base2-form}) over $\C$. {}From these ingredients, 
we build $\extd_1:\Omega^1(H)\to\Omega^2(H)$ as a $96\times 48$-matrix 
defined by $\extd_1(fe_a)=\extd_0(f)\wedge e_a+f\extd e_a$. We then compute 
the kernel of $\extd_1$ and find it to be 12-dimensional. The image of 
$\extd_0$ is necessarily 11-dimensional (its kernel is the constant 
functions) and hence the cohomology is 1-dimensional.
$\theta$ is closed but never exact (for any finite group) and hence 
represents this class. \eproof

Next, whereas the cohomology is a linear problem, one can also consider its 
non-linear
variant called $U(1)$-gauge theory. Here we define the curvature of a 1-form 
$\alpha\in \Omega^1(H)$ to be the 2-form $F(\alpha)=\extd\alpha+\alpha\wedge
\alpha$. This transforms by conjugation under the gauge transform 
$\alpha\mapsto u\alpha u^{-1}+u\extd u^{-1}$ for any no-where zero function 
$u\in H$. One can also impose here
unitarity conditions as in \cite{mr}. In this context it would be 
interesting to
find the moduli space of (unitary) flat connections. This was done for $S_3$
in \cite{mr} and found to have a richer structure than the cohomology alone, 
i.e. with other solutions beyond multiples of $\theta$ and we would expect 
something similarly rich for $A_4$. For example, if we focus on flat 
connections with constant coefficients in the 
$\{e_a\}$ basis as in \cite{maj-per}, a short computation shows that these are 
given by the five lines
\eqn{flat}{ \lambda e_t-\theta,\quad \lambda e_x-\theta, 
\quad\lambda e_y-\theta,\quad 
\lambda e_z-\theta,\quad (\lambda-1)\theta}
for $\lambda$ a parameter. There is a similar behaviour for any cyclic 
conjugacy class.

A further question relates to the fact that $A_4\subset S_4$ as a normal 
subgroup. Therefore its exterior algebra should be related to that of $S_4$ 
for a suitable conjugacy class on that. The different differential structures 
on $S_4,S_5$ for different conjugacy classes are studied in \cite{maj-per} 
and looking there, one finds that the order 8 conjugacy class containing 
$(123)$ in $S_4$ has the required exterior algebra. Its eight generators 
split into two 
sets, namely $\{e_{(123)},e_{(134)},e_{(243}),e_{(142)}\}$ generating a 
subalgebra with relations as in Proposition~\ref{prop3.1} and a complementary 
set $\{e_{(132)},e_{(143)},e_{(234)},e_{(124)}\}$ generating the opposite 
subalgebra. We denote the latter generators by $\{e_{\bar t}, e_{\bar x},
e_{\bar y}, 
e_{\bar z}\}$. There are nontrivial cross relations between the two sets: 
\eqn{crossext}{ e_a\wedge e_{\bar a}+e_{\bar a}\wedge e_a=0,}
\[ e_t\wedge e_{\bar z}+e_{\bar z}\wedge e_x+e_x\wedge e_{\bar y}+e_{\bar y}
\wedge e_t=0\]
\[ e_t \wedge e_{\bar x}+e_{\bar x}\wedge e_y+e_y\wedge e_{\bar z}
+e_{\bar z}\wedge e_t=0\]
\[ e_t\wedge e_{\bar y}+e_{\bar y}\wedge e_z+e_z\wedge e_{\bar x}+e_{\bar x}
\wedge e_t=0\]
and their three conjugates (given by applying $\bar{\ }$ and reversing 
products). 
In other words, the differential geometry of $S_4$ appears to be some form 
of `complexification' of that of $A_4$. 

Indeed, the conjugacy class $\{\bar t, \bar x, \bar y, \bar z\}$ in $A_4$ 
is also cyclic and the equations of its associated
exterior algebra are of the same form as in Proposition~\ref{prop3.1}. 
It defines a differential geometry on $A_4$ conjugate to the one we have 
studied above. Indeed, one has
\eqn{conjdel}{ \del_a^\dagger=\del_{a^{-1}}=\del_{a^2}=\del_{\bar a}}
where ${}^\dagger$ denotes transpose with respect to the $l^2$ inner 
product on $A_4$. Here the first equalitity is a general
feature for any finite group and follows from the braided-Leibniz rule. The 
second equality is due to all elements of $\CC$ in our case having order 3 
and the third is a special 
feature $a^2=\bar a$ of 
the multiplication Table~3. To complete the analysis 
let us note that $A_4$ has just one other nontrivial conjugacy class, 
$\{u,v,w\}$. This does not generate $A_4$ i.e., the quantum manifold 
structure that it defines is not connected (not every point can be reached 
from any other by steps taken from the conjugacy 
class). The component of this connected to the group identity is 
$\Z_2\times \Z_2$ with its universal differential calculus.

Finally, returning to the Riemannian geometry, one can and should consider 
more general metrics and vierbeins. Since we have found above a canonical 
invariant $\eta$, one could use this to fix the relationship between 
a general vierbein and covierbein, i.e. look for metrics of the form 
$g=\eta^{ab}e_a\tens_H e_b$ with $\{e_a\}$ as the free variable (rather 
than $e,e^*$ independent). Among the moduli space of pairs $(e,A)$ of vierbein 
and spin connection (or more generally of triples $(e,e^*,A)$), our results 
above show that there is at least one canonical point where the Ricci tensor
vanishes. This motivates the problem of solving the Ricci flat equations in 
general, i.e. classical `gravity' on $A_4$. Similarly, one can consider
functional integrals, now a finite number of usual integrals, with 
Einstein-Hilbert action, i.e. `quantum-gravity'. These are difficult 
nonlinear questions unsolved even for
$S_3$ and beyond our present scope.

\subsection*{Acknowledgements} One of us (F.N) would like to thank 
 Prof. J-P Antoine for facilitating his stay at FYMA-UCL during year
2000-2001.


\begin{thebibliography}{15}
\bibitem{maj-br}
S. Majid. Quantum and braided group Riemannian
 geometry. J. Geom. Phys., {\bf 30}:113-146,1999.

\bibitem{maj}
 S. Majid. Riemannian geometry on quantum groups and finite groups
 with nonuniversal differentials. preprint, June 2000.

\bibitem{w}
S. L. Woronowicz. Differential Calculus on Compact Matrix
Pseudogroups (Quantum Groups) Commun. Math. Phys. {\bf 122}:
125-170 (1989).

\bibitem{tbmaj} T. Brzezi\'nski and S. Majid.
Quantum group gauge theory on quantum spaces. Commun. Maths.
Phys., {\bf 157}:591-638, 1993. Erratum {\bf 167}:235, 1995.


\bibitem{mr}
 S. Majid and E. Raineri. Electromagnetism and gauge theory on the
 permutation group $S_{3}$. Preprint, Dec. 2000.

\bibitem{con}
A. Connes. Noncommutative Geometry, Academic Press.

\bibitem{maj-per}
 S. Majid. Noncommutative differentials and Yang-Mills on
 permutation groups $S_N$. Preprint, May 2001.



\bibitem{dim}
A. Dimakis and F. Muller-Hoissen. Discrete Riemannian geometry. J.
Math. Phys., {\bf 40}: 1518, 1999.

\bibitem{kbres}
K. Bresser, F. Muller-Hoissen, A. Dimakis and A. Sitarz.
Noncommutative geometry of finite groups. J. Phys. A. Math. Gen.,
{\bf 29}: 2705-2735, 1996.

\bibitem{sitarz}
A. Sitarz. Noncommutative geometry and gauge theory on discrete
groups. J. Geom. Phys., {\bf 15}: 123-136,1995.

\bibitem{maj-con}
S. Majid. Conceptual issues for noncommutative gravity on algebras
and finite sets, Int. J. Mod. Phys.{\bf B 14} (2000) 2427--2449.

\bibitem{maj-book}
S. Majid. Foundations of Quantum Group Theory, Cambridge Univ. Press, 1995.

\end{thebibliography}
\end{document}